\documentclass[12pt]{article}
\usepackage{fancyhdr, array,calc,graphicx,url,tabularx}
\usepackage{geometry}
\usepackage{latexsym}
\usepackage{amssymb}
\usepackage{amsmath}
\usepackage{makeidx}
\usepackage{enumerate}
\usepackage{verbatim}
\usepackage{tikz}
\usetikzlibrary{calc,arrows.meta,positioning}
\usepackage[colorlinks=true,linkcolor=blue]{hyperref}
\usepackage{changepage}

\makeindex

{
   \end{minipage}
   \vspace*{\stretch{3}}
   \clearpage
}

\newcommand{\bR}{{\mathbb{R}}}

\newcommand{\rest}{\restriction}

\newcommand{\card}[1]{{\vert #1 \vert} }

\newcommand{\forces}{\Vdash}

\renewcommand{\models}{\vDash}
\newcommand{\powerset}{{\wp}}
\newcommand{\cp}{{\rm crit }}

\newtheorem{theorem}{Theorem}[section]

\newtheorem{definition}[theorem]{Definition}

\newtheorem{lemma}[theorem]{Lemma}
\newtheorem{corollary}[theorem]{Corollary}
\newtheorem{claim}[theorem]{Claim}

\newtheorem{remark}[theorem]{Remark}

\numberwithin{figure}{section}

\newenvironment{proof}{{\it{
Proof.}}}{\nopagebreak\mbox{}{\hfill$\square$}
\par\bigskip}
\newenvironment{Claim}{\noindent{\it
Claim} {}}{}
\newenvironment{Case}{\noindent{\it
Case}}{}

\newcommand{\rthm}[1]{Theorem~\ref{#1}}
\newcommand{\rlem}[1]{Lemma~\ref{#1}}

\newcommand{\rcor}[1]{Corollary~\ref{#1}}
\newcommand{\rdef}[1]{Definition~\ref{#1}}

\newcommand{\rsec}[1]{Section~\ref{#1}}

\def\k{\kappa}
\def\a{\alpha}
\def\b{\beta}
\def\d{\delta}

\def\l{\lambda}

\def\P{{\mathcal{P} }}
\def\W{{\mathcal{W} }}
\def\Q{{\mathcal{ Q}}}

\def\R{{\mathcal R}}
\def\X{{\mathbb X}}

\def\M{{\mathcal{M}}}
\def\N{{\mathcal{N}}}

\def\T {{\mathcal{T}}}
\def\U{{\mathcal{U}}}
\def\S{{\mathcal{S}}}
\def\V{{\mathcal{V}}}
\def\X{{\mathcal{X}}}

\def\card#1{\left|#1\right|}

\def\iff{\mathrel{\leftrightarrow}}

\def\and{\mathrel{\kern1pt\&\kern1pt}}

\def\<#1>{\langle\,#1\,\rangle}

 \input xy
 \xyoption{all}
       
\newcounter{nameOfYourChoice}

\title{$\sf{Sealing}$ from iterability
}

\author{Grigor Sargsyan \footnote{Department of Mathematics, Rutgers University, NJ, USA. Email: gs481@math.rutgers.edu} 
\\ Nam Trang \footnote{Department of Mathematics, University of North Texas, Denton, TX, USA. Email: Nam.Trang@unt.edu}}

\date{October 12, 2019\footnote{Revised December 15, 2020}}

\begin{document}

\maketitle

\begin{abstract}
We show that if $V$ has a proper class of Woodin cardinals, a strong cardinal, and a generically universally Baire iteration strategy (as defined in the paper) then $\sf{Sealing}$ holds after collapsing the successor of the least strong cardinal to be countable. This result is complementary to \cite[Theorem 3.1]{StrengthuB} where it is shown that $\sf{Sealing}$ holds in a generic extension of a certain minimal universe. The current theorem is more general in that no minimality assumption is needed. A corollary of the main theorem is that $\sf{Sealing}$ is consistent relative to the existence of a Woodin cardinal which is a limit of Woodin cardinals. This improves significantly on the first consistency of $\sf{Sealing}$ obtained by W.H. Woodin. 

The $\sf{Largest\ Suslin\ Axiom}$ ($\sf{LSA}$) is a determinacy axiom isolated by Woodin. It asserts that the largest Suslin cardinal is inaccessible for ordinal definable bijections. Let $\sf{LSA-over-uB}$ be the statement that  in all (set) generic extensions there is a model of $\sf{LSA}$ whose Suslin, co-Suslin sets are the universally Baire sets. The other main result of the paper shows that assuming $V$ has a proper class of inaccessible cardinals which are limit of Woodin cardinals, a strong cardinal, and a generically universally Baire iteration strategy, in the universe $V[g]$, where $g$ is $V$-generic for the collapse of the successor of the least strong cardinal to be countable, the theory $\sf{LSA-over-UB}$ fails; this implies that $\sf{LSA-over-UB}$ is not equivalent to $\sf{Sealing}$ (over the base theory of $V[g]$).\footnote{This is interesting and somewhat unexpected, in light of \cite[Theorem 1.6]{StrengthuB}. Compare this result with Steel's well-known theorem that ``$\sf{AD}^{L(\mathbb{R})}$ holds in all generic extensions" is equivalent to ``the theory of $L(\mathbb{R})$ is sealed" in the presence of a proper class of measurable cardinals.}
\end{abstract}

We identify elements of the Baire space $\omega^\omega$ with reals. Throughout the paper, by a ``set of reals $A$", we mean $A\subseteq \omega^\omega$. A set of reals $A$ is \textit{$\gamma$-universally Baire} if there are trees $T,U$ on $\omega\times \lambda$ for some $\lambda$ such that $A = p[T] = \mathbb{R}\backslash p[U]$ and whenever $g$ is a $< \gamma$-generic, in $V[g]$, $p[T] = \mathbb{R}\backslash p[U]$. We write $A^g$ for $p[T]^{V[g]}$; this is the canonical interpretation of $A$ in $V[g]$.\footnote{One can show $A^g$ does not depend on the choice of $T,U$.} $A$ is \textit{universally Baire} if $A$ is $\gamma$-universally Baire for all $\gamma$. Let $\Gamma^\infty$ be the set of universally Baire sets. Given a generic $g$, we let $\Gamma^\infty_g=(\Gamma^\infty)^{V[g]}$ and $\bR_g=\bR^{V[g]}$. The next definition is due to Woodin.

\begin{definition}\label{dfn:ub_sealing} $\sf{Sealing}$ is the conjunction of the following statements.
\begin{enumerate}
\item For every set generic $g$, $L(\Gamma^\infty_g, \mathbb{R}_g)\models \sf{AD}^+$ and $\powerset(\bR_g)\cap L(\Gamma^\infty_g, \mathbb{R}_g)=\Gamma^\infty_g$.
\item  For every set generic $g$ over $V$, for every set generic $h$ over $V[g]$, there is an elementary embedding 
\begin{center}
$j: L(\Gamma^\infty_g, \mathbb{R}_g)\rightarrow L(\Gamma^\infty_h, \mathbb{R}_h)$.
\end{center}
\end{enumerate}
 such that for every $A\in \Gamma^\infty_g$, $j(A)=A^h.$  
 \end{definition}
 
 $\sf{Sealing}$ is a form of Shoenfield-type generic absoluteness for the theory of universally Baire sets. In this paper, we will avoid motivational discussion as \cite{StrengthuB} has a lengthy introduction to the subject. We should say, however, that $\sf{Sealing}$ is an important hypothesis in set theory and particularly in inner model theory for several reasons. If a large cardinal theory $\phi$ implies $\sf{Sealing}$ then the Inner Model Program for building canonical inner models of $\phi$ cannot succeed (at least with the criteria for defining ``canonical inner models" as is done to date), cf \cite[Sealing Dichotomy]{StrengthuB}. $\sf{Sealing}$ signifies a place beyond which new methodologies are needed in order to advance the Core Model Induction techniques. In particular, to obtain consistency strength beyond $\sf{Sealing}$ from strong theories such as the Proper Forcing Axiom, one needs to construct canonical subsets of $\Gamma^\infty$ (third-order objects), instead of elements of $\Gamma^\infty$ like what has been done before (see \cite[Section 1]{StrengthuB} for a more detailed discussion). The consistency of $\sf{Sealing}$ was first demonstrated by Woodin, who showed that if there is a proper class of Woodin cardinals and a supercompact cardinal $\k$ then $\sf{Sealing}$ holds after collapsing $2^{2^\k}$ to be countable. Woodin's proof can be found in \cite{StationaryTower}. 
 
One of the main corollaries of the \rthm{main theorem} is that the set theoretic strength of $\sf{Sealing}$ is below a Woodin cardinal that is a limit of Woodin cardinals; this improves significantly the aforementioned result of Woodin. Another proof of this fact was presented in \cite{StrengthuB}, where the authors establish an actual equiconsistency for $\sf{Sealing}$. One advantage of the proof in this paper is that no smallness assumption is made (unlike \cite{StrengthuB}). Another, perhaps more important, advantage of the current proof over the one presented in \cite{StrengthuB} is that this proof is more accessible. Our proof of $\sf{Sealing}$ is based on iterability and uses recent ideas from descriptive inner model theory. However, in this paper, our aim is to present the proof of our main theorem, \rthm{main theorem}, without using any fine structure theory or heavy machinery from inner model theory, so that the paper is accessible to the widest possible audience.  We will only assume general knowledge of iterations, iteration strategies and Woodin's extender algebra, all of which are topics that can be presented without any fine structure theory. For instance, the reader can consult \cite{FarahEA} or \cite{IT}. The fact that the hypothesis of \rthm{main theorem} is weaker than a Woodin cardinal that is a limit of Woodin cardinals follows from a very recent work of Steel (\cite{normalization_comparison}) and the first author (\cite{GG}) (but also see \cite{RecentResultsIMT}), and this fact will not be proven here, as it is well beyond the scope of this paper.

Given a transitive model $Q$ of set theory and a $Q$-cardinal $\k$, we let $Q|\k=H_\k^Q$. We say $E$ is a $(\k, \l)$-short extender over $Q$ if there is a $\Sigma_1$-elementary embedding $j:Q\rightarrow M$ such that 
\begin{enumerate}
\item $M$ is transitive,
\item $M=\{j(f)(a): f: \k^{<\omega} \rightarrow Q, f\in Q$ and $a\in \l^{<\omega}\}$,
\item $j(\k)\geq \l$, and
\item $E=\{ (a, A): a\in \l^{<\omega}, A\subseteq [\k]^{\card{a}}$ and $a\in j(A)\}$.
\end{enumerate}
$\kappa$ is called the \textit{critical point of $E$} and $\lambda$ the \textit{length of $E$}. We write $\kappa = \cp(E)$ and $\lambda = lh(E)$. $M$ is then called the \textit{ultrapower of $Q$ by $E$} and is uniquely determined by $Q$ and $E$. We write $M = Ult(Q,E)$. Given a set $X$ and an extender $E$, we say $E$ \textit{coheres} $X$ if $X\cap V_{lh(E)}=j(X)\cap V_{lh(E)}$. For more on short extenders, the reader can consult \cite{IT} or \cite{FarahEA}.

We can also define the notion of \textit{a long extender}, though we will not need the precise definition in this paper. Roughly speaking, given an elementary embedding $j:V\rightarrow M$ with critical point $\kappa$, an ordinal $\eta>\kappa$, and letting $\xi$ be least such that $j(\xi)\geq \eta$, we can define an extender $E$ of length $\eta$ from $j$. This is a function $F:\powerset(\xi)\rightarrow V$ given by: $F(A) = j(A)\cap \eta$. If $\xi > \kappa$, then $E$ is a long extender. For more details on long extenders, see \cite{woodin2010suitable}.

Suppose $P$ is a transitive model of set theory. We let ${\sf{ile}}(P)$ be the set of \textit{inaccessible-length extenders} of $P$. More precisely $\sf{ile}(P)$ consists of short extenders $E\in P$ such that $P\models ``lh(E)$ is inaccessible and $V_{lh(E)}=V_{lh(E)}^{Ult(V, E)}$."

\begin{definition}\label{pre-iterable} We say that $\P$ is a pre-iterable structure if $\P=(P, {\sf{ile}}(P))$ where $P$ is a transitive model of $\sf{ZFC}$. 
\end{definition}

 When we talk about iterability for $\P$, we mean iterability with respect to extenders in $\vec{E}^\P=_{def}{\sf{ile}}(P)$ (and its images). Thus, the relevant iterations are those that are built by using extenders in $\vec{E}$ and its images.

Recall from \cite{IT} that an iteration $\T$ is \textit{normal} if the extenders used in it have increasing lengths and each extender $E$ used along $\T$ is applied to the least possible model, i.e. $E$ is applied to the first model $\M^\T_\alpha$ where the ultrapower $Ult(\M^\T_\alpha,E)$ makes sense. Following Jensen, we will say that $\T$ is a \textit{smooth iteration} (of its base model) if it can be represented as a \textit{stack} of normal iterations. More precisely, $\T=(\T_i: i<\eta)$ where $\T_0$ is a normal iteration of the base model of $\P$ and for $i\in (0, \eta)$, $\T_i$ is a normal iteration of the last model of $\T_{i-1}$ if $i$ is a successor ordinal and on the direct limit of $(\T_j : j<i)$ under the iteration embeddings if $i$ is limit. We say that a pre-iterable structure $\P$ is \textit{smoothly iterable} if player II has a wining strategy in the iteration game of arbitrary length that produces smooth iterations. Recall that in iteration games, player I picks the extenders while player II plays branches at limit steps. We say that $\Sigma$ is an iteration strategy for $\P$ if it is a strategy for $\P$ in the iteration game that produces arbitrary length smooth iterations of $\P$. 

Finally we state self-iterability. The Unique Branch Hypothesis ($\textsf{UBH}$) is the statement that every normal iteration tree $\T$ on $V$ has at most one cofinal well-founded branch. The Generic Unique Branch Hypothesis ($\sf{gUBH}$) says that $\sf{UBH}$ holds in all set generic extensions. The notion of generically universally Baire (guB) strategy appears in the next section as \rdef{ub strategy}.

\begin{definition}\label{self-iterability} We say that self-iterability holds if the following holds in $V$.
\begin{enumerate}
\item $\sf{gUBH}$.
\item $\V=(V, \sf{ile}(V))$ is a pre-iterable structure that has a guB-iteration strategy.
\end{enumerate}
\end{definition}
Notice that because of clause 1, the iteration strategy in clause 2 is unique.

\begin{theorem}\label{main theorem} Assume self-iterability holds, and suppose there is a class of Woodin cardinals and a strong cardinal. Let $\k$ be the least strong cardinal of $V$ and let $g\subseteq Coll(\omega, \k^+)$ be $V$-generic. Then $V[g]\models \sf{Sealing}$.
 \end{theorem}
 
As mentioned above, a corollary of \rthm{main theorem}, via a non-trivial amount of work in \cite{normalization_comparison} and \cite{GG} (but also see \cite{RecentResultsIMT}\footnote{The existence of an lbr hod premouse $\P$ as in \cite[Theorem 1.2]{RecentResultsIMT} follows from the existence of a Woodin limit of Woodin cardinals by \cite[Step 4]{RecentResultsIMT}. Then letting $\lambda_0$ be as in \cite[Theorem 1.2]{RecentResultsIMT}, $\P|\lambda_0$ satisfies the hypothesis of Theorem \ref{main theorem}.}), is
\begin{corollary}\label{cor:wlw}
Con$(\sf{ZFC}\ $+ there is a Woodin cardinal which is a limit of Woodin cardinals$)$ implies Con$(\sf{Sealing})$.
\end{corollary} 
 
The main idea behind the proof of \rthm{main theorem} originates in \cite{StrengthuB}. The most relevant portion of that paper is \cite[Theorem 3.1]{StrengthuB}. We should note that the hypothesis of \rthm{main theorem} cannot be weakened to just $\sf{gUBH}$ for plus-2 iterations as this form of $\sf{UBH}$ holds in a minimal mouse with a strong cardinal, a class of Woodin cardinals and a stationary class of measurable cardinals\footnote{This fact is due to Steel, see \cite[Theorem 3.3]{LKC}.}, but this theory is weaker than $\sf{Sealing}$ as shown  by \cite[Theorem 3.1]{StrengthuB}.

The $\sf{Largest\ Suslin\ Axiom}$ was introduced by Woodin in \cite[Remark 9.28]{Woodin}. The terminology is due to the first author. Here is the definition. In the following, we say that a cardinal $\kappa$ is \textit{$\sf{OD}$-inaccessible} if for every $\alpha<\kappa$ there is no surjection $f: \powerset(\alpha)\rightarrow \kappa$ that is definable from ordinal parameters.
\begin{definition}\label{dfn:lsa} 
The $\sf{Largest\ Suslin\ Axiom}$, abbreviated as $\sf{LSA}$, is the conjunction of the following statements:
\begin{enumerate}
\item $\sf{AD}^+$.
\item There is a largest Suslin cardinal.
\item The largest Suslin cardinal is $\sf{OD}$-inaccessible.
\end{enumerate}
 \end{definition}
 
 In the hierarchy of determinacy axioms, which one may appropriately call the $\sf{Solovay\ Hierarchy}$\footnote{Solovay defined what is now called the $\sf{Solovay\ Sequence}$ (see \cite[Definition 9.23]{Woodin}). It is a closed sequence of ordinals with the largest element $\Theta$, where $\Theta$ is the least ordinal that is not a surjective image of the reals. One then obtains a hierarchy of axioms by requiring that the $\sf{Solovay\ Sequence}$ has complex patterns. $\sf{LSA}$ is an axiom in this hierarchy. The reader may consult \cite{BSL} or \cite[Remark 9.28]{Woodin}.}, $\sf{LSA}$ is an anomaly as it belongs to the successor stage of the $\sf{Solovay\ Hierarchy}$ but does not conform to the general norms of the successor stages of the $\sf{Solovay\ Hierarchy}$. Prior to \cite{hod_mice_LSA}, $\sf{LSA}$ was not known to be consistent. \cite{hod_mice_LSA} shows that it is consistent relative to a Woodin cardinal that is a limit of Woodin cardinals. Nowadays, the axiom plays a key role in many aspects of inner model theory, and features prominently  in Woodin's $\sf{Ultimate\ L}$ framework (see \cite[Definition 7.14]{UltimateL} and Axiom I and Axiom II on page 97 of \cite{UltimateL}\footnote{The requirement in these axioms that there is a strong cardinal which is a limit of Woodin cardinals is only possible if $L(A, \bR)\models \sf{LSA}$.}). 
 
  \begin{definition}\label{dfn:generic_LSA}
Let $\sf{LSA-over-uB}$ be the statement: For all $V$-generic $g$, in $V[g]$, there is $A\subseteq \bR_g$ such that $L(A, \bR_g)\models \sf{LSA}$ and $\Gamma^\infty_g$ is the Suslin co-Suslin sets of $L(A, \bR_g)$.  
\end{definition}

\cite{StrengthuB} shows that $\sf{Sealing}$ is equiconsistent with $\sf{LSA-over-UB}$ over the theory $\sf{ZFC} +\ $``there is a proper class of Woodin cardinals and the class of measurable cardinals is stationary". In this paper, we show that in general, one cannot replace ``equiconsistent" with ``equivalent". Recall from \cite{normalization_comparison} the statement of \textit{Hod Pair Capturing} $(\sf{HPC})$: for any Suslin co-Suslin set $A$, there is a least-branch (lbr) hod pair $(\P,\Sigma)$ such that $A$ is definable from parameters over $(HC,\in, \Sigma)$. \textit{No Long Extender} $(\sf{NLE})$ is the statement: there is no countable, $\omega_1+1$-iterable pure extender premouse $M$ such that there is a long extender on the $M$-sequence. The notion of least-branch hod mice (lbr hod mice) is defined precisely in \cite[Section 5]{normalization_comparison}.

\begin{definition}\label{dfn:gHPC}
$\sf{gHPC}$ is the statement: suppose $V[g]$ is a set generic extension of $V$, suppose in $V[g]$, $M = L(\Gamma,\mathbb{R})$ is a model of $\sf{AD}^+$. Then $M\vDash \sf{HPC}$.
\end{definition}

\begin{theorem}\label{thm:not_equiv}
Suppose self-iterability holds and there is a proper class of inaccessible cardinals which are limit of Woodin cardinals. Suppose $\sf{gHPC}$ and $\sf{NLE}$ hold. Then $V\models \sf{LSA-over-UB}$ fails.
\end{theorem}

\begin{remark}\label{rmk:Wdn}
\begin{enumerate}
\item The hypotheses of Theorem \ref{thm:not_equiv} hold in the universe of lbr hod mice that have a proper class of inaccessible cardinals which are limit of Woodin cardinals (cf. \cite{normalization_comparison}). So such hod mice satisfy ``$\sf{LSA-over-UB}$ fails."
\item Woodin has independently shown that $\sf{LSA-over-UB}$ can fail. More precisely, $\sf{LSA-over-UB}$ fails assuming there is a proper class of Woodin cardinals, a proper class of strong cardinals, and there is an inaccessible cardinal which is a limit of Woodin cardinals and strong cardinals.
\end{enumerate}
\end{remark}

Remark \ref{rmk:Wdn}(1), Theorem \ref{thm:not_equiv}, and the fact that self-iterability and $\sf{gHPC}$ hold in any generic extension of an lbr hod mouse with a proper class of Woodin cardinals give us the following.
\begin{corollary}\label{cor:not_equiv}
Let $V$ be the universe of an lbr hod mouse with a proper class of inaccessible cardinals which are limit of Woodin cardinals, and a strong cardinal.  Assume $\sf{NLE}$. Let $\kappa$ be the least strong cardinal of $V$ and $g\subseteq Coll(\omega,\kappa^+)$ be $V$-generic. Then $V[g] \models \sf{Sealing}$ holds and $\sf{LSA-over-UB}$ fails.
\end{corollary}

\rcor{cor:not_equiv} is surprising. For example, generic absoluteness for $L(\bR)$, namely that for all successive generics $g$ and $h$ there is an elementary embedding $j:L(\bR_g)\rightarrow L(\bR_{g*h})$, is equivalent to the existence and the universally Bairness of the next canonical set beyond $L(\bR)$, namely $\bR^\#$\footnote{This fact is due to Steel and Woodin. For example, see genericity iterations in \cite{steel2010outline}.}. While one cannot hope that $\sf{Sealing}$ would imply both the existence and the universal Bairness of the next canonical set of reals beyond $\Gamma^\infty$\footnote{As all universally Baire sets are already in $\Gamma^\infty$.}, one could still hope that the cause of $\sf{Sealing}$ is the existence of some nice set of reals just like the cause of the generic absoluteness of $L(\bR)$ is the universally Bairness of $\bR^\#$\footnote{Or rather the universal Bairness of the $\omega_1$-iteration strategy of $\M_{\omega}^\#$.}. Because the next nice set beyond $\Gamma^\infty$ cannot be universally Baire, the best we can hope for is that the next set beyond $\Gamma^\infty$ creates an $\sf{LSA}$ model over $\Gamma^\infty$. In fact, this discussion was the original motivation for isolating $\sf{LSA-over-UB}$. However, unlike our expectations, what causes $\sf{Sealing}$ may not be coded into a set of reals as demonstrated by \rcor{cor:not_equiv}.  

Throughout this paper, except in \rsec{sec:gub},  we assume the hypothesis of \rthm{main theorem}. Throughout this paper, except in \rsec{sec:gub}, $\k$ will stand for the least strong cardinal. In this paper, especially in \rsec{sec: gen-it}, we will make heavy use of Neeman's ``realizable maps are generic" result that appears as  \cite[Corollary 4.9.2]{Neeman}. Sections \ref{sec: cap ub} and \ref{sec:der model} make heavy use of the results of Section \ref{sec: gen-it} to show that for $V$-generic $g\subseteq Coll(\omega,\kappa^+)$, where $\kappa$ is as in \rthm{main theorem}, for $V[g]$ generic $h$, one can realize $L(\Gamma_{g*h}^\infty,\mathbb{R}_{g*h})$ as the derived model of an iterate of a countable substructure of $V_\gamma[g*h]$ for some large $\gamma$ (Lemma \ref{der model rep}). This is then used to prove \rthm{main theorem} in Section \ref{sec:main}. The last section proves Theorem \ref{thm:not_equiv}.

\textbf{Acknowledgement.} The authors would like to thank the NSF for its generous support. The first author is supported by NSF Career Award DMS-1352034. The second author is supported by NSF Grants DMS-1565808, DMS-1849295, and NSF Career Grant DMS-1945592.

\section{Generically universally Baire iteration strategies}\label{sec:gub}

In this paper we will need three properties of iteration strategies, namely \textit{Skolem-hull condensation}, \textit{pullback condensation} and \textit{generically universal Bairness}. We now define these notions. 

We say $(\P, \Psi)$ is an \textit{iterable pair} if $\P$ is a pre-iterable structure and $\Psi$ is a strategy for it. Suppose $(\P, \Psi)$ is an iterable pair. If $\T$ is a smooth iteration of $\P$ according to $\Psi$ with last model $\Q$ then we write $\Psi_{\T, \Q}$ for the strategy of $\Q$ induced by $\Psi$. Namely, $\Psi_{\T, \Q}(\U)=\Psi(\T^\frown \U)$. When $\Psi_{\T, \Q}$ is independent of $\T$ we will drop it from our notation. Given a $\P$-cardinal $\xi$, we write $\Psi_{\P|\xi}$ for the fragment of $\Psi$ that acts on smooth iterations based on $\P|\xi$. Here recall that $\P|\xi=H_\xi^\P$.  

Continuing with $(\P, \Psi)$, suppose $\pi:\N\rightarrow \P$ is elementary. Given a smooth iteration $\T$ of $\N$ we can define the copy $\pi\T$ on $\P$ which may or may not have well-founded models. The construction of $\pi\T$ was introduced in \cite{IT} on page 17. Suppose now that $\T$ is such that $\pi\T$ is according to $\Psi$ and $\T$ is of limit length. Let $b=\Psi(\pi\T)$. It follows from the construction of $\pi\T$ that $b$ yields a well-founded branch of $\T$. 

We then say $\Lambda$ is the $\pi$-pullback of $\Psi$ if for any smooth iteration $\T$ on $\N$ that is according to $\Lambda$, $\pi\T$ is according to $\Psi$. It is customary to let $\Lambda$ be $\Psi^\pi$. 

\begin{definition}\label{skolem hull condensation} Suppose $(\P, \Psi)$ is an iterable pair. 
 We say $\Psi$ has \textbf{Skolem-hull condensation} if whenever $\T$ is an iteration according to $\Psi$, $\xi$ is such that $\T\in V_\xi$ and $\pi: M\rightarrow V_\xi$ is elementary such that $(\P|\xi, \Psi_{\P|\xi}, \T)\in rng(\pi)$ then $\pi^{-1}(\T)$ is according to $\Psi^\pi_{\P|\xi}$. 
 \end{definition}
 
 \begin{definition}\label{pullback condensation} Suppose $(\P, \Psi)$ is an iterable pair. We say $\Psi$ has \textbf{pullback condensation} if whenever $\T$ is an iteration according to $\Psi$ with last model $\Q$ and $\U$ is an iteration of $\Q$ according to $\Psi_{\T, \Q}$ with last model $\R$ then $\Psi^{\pi^\U}_{\T^\frown \U, \R}=\Psi_{\T, \Q}$.
  \end{definition}
  
  The following theorems are easy consequences of $\sf{UBH}$ ($\sf{gUBH}$), and are probably not due to the authors. 
  \begin{theorem}\label{easy consequence0} Assume $\sf{UBH}$ and suppose $\l$ is inaccessible. Then $V_\l\models \sf{UBH}$. 
  \end{theorem}
  \begin{theorem}\label{easy consequence1} Assume self-iterability and suppose $\Psi$ is the unique strategy of $\V$. Then $\Psi$ has Skolem-hull condensation and pullback condensation. 
  \end{theorem}

Suppose $(\P, \Psi)$ is an iterable pair. Given a strong limit cardinal $\k$ and $F\subseteq Ord$, set 
\begin{center}
$W^{\Psi, F}_\k=(H_{\k}, F\cap \kappa, \P|\k, \Psi_{\P|\k}\rest H_\k, \in)$.
\end{center}

Given a structure $Q$ in a language extending the language of set theory with a transitive universe, and an $X\prec Q$, we let $ M_X$ be the transitive collapse of $X$ and $\pi_X: M_X\rightarrow Q$ be the inverse of the transitive collapse. In general, the preimages of objects in $X$ will be denoted by using $X$ as a subscript, e.g. $\pi_X^{-1}(\P) = \P_X$. Suppose in addition $Q=(R,...\P,\Phi,...)$ where $\P$ is a pre-iterable structure and $\Phi$ is an iteration strategy of $\P$. We will then write $X\prec (Q|\Phi)$ to mean that $X\prec Q$ and the strategy of $\P_X$ that we are interested in is $\Phi^{\pi_X}$. We set $\Lambda_X=\Phi^{\pi_X}$.

Motivated by the definition of universally Baire sets that involves club of generically correct hulls, we make the following definition.

\begin{definition}\label{ub strategy} We say $\Psi$ is a \textbf{generically universally Baire (guB) strategy} for a pre-iterable $\P=(P, \vec{E})$ if there is a formula $\phi(x)$ in the language of set theory augmented by three relation symbols and $F\subseteq Ord$ such that for every inaccessible cardinal $\k$ and for every countable
\begin{center}
$X\prec (W^{\Psi, F}_\k| \Psi_{\P|\k})$ 
\end{center}
 whenever 
\begin{enumerate}[(a)]
\item $g\in V$ is $M_{X}$-generic for a poset of size $<\k_X$ and 
\item $\T\in M_X[g]$ is such that for some $M_X$-inaccessible $\eta<\k_X$, $\T$ is an iteration of $\P_X|\eta$,
\end{enumerate}
 the following conditions hold:
\begin{enumerate}
\item if $lh(\T)$ is a limit ordinal and $\T\in dom(\Lambda_X)$ then $\Lambda_X(\T)\in M_X[g]$,
\item $\T$ is according to $\Lambda_X$ if and only if $\M_X[g]\models \phi[\T]$.
\end{enumerate}
We say that $(\phi, F)$ is a generic prescription of $\Psi$.
\end{definition}
In \rdef{ub strategy}, we could demand that there is a club of $X$ with the desired properties. However that would be equivalent to our definition as we can let $F$ above code the desired club. In the next section our goal is to prove some basic facts about $guB$-strategies.

\section{Generic interpretability of guB strategies}\label{sec: gen-it}

As we said in the introduction, from this point on we work under the hypothesis of \rthm{main theorem}. However, we will not use the existence of a strong cardinal until \rsec{sec:der model}.

 Let $\Psi$ be the guB-strategy of $\V=(V, \sf{ile}(V))$ and fix a generic prescription $(\phi, F)$ for $\Psi$ (see Definition \ref{ub strategy}). We will omit $\Psi, F$ from our notation and just write $W_\k$ instead of $W_\k^{\Psi, F}$. Given a cardinal $\a$ we will write $\Psi_\a$ for the fragment of $\Psi$ that acts on iterations based on $\V|\a$. Often we will treat $\Psi_\a$ as a strategy for $\V|\a$ rather than a strategy for $\V$. Similarly, given an interval $(\a, \b)$ we will write $\Psi_{\a, \b}$ for the fragment of $\Psi$ on iterations based on $\V|\beta$ above $\alpha$. To make the notation simpler, often we will not specify the domain of $\Psi_\a$ that we have in mind (as in \rlem{simple capturing}). 

Let $\d$ be a Woodin cardinal of $\V$. We first prove that $\Psi_\d$ has canonical extensions in generic extensions of $V$. As a first step, we prove the following useful capturing result.

\begin{lemma}\label{simple capturing} Suppose $\l$ is an inaccessible cardinal and let $X\prec (W_\l|\Psi_\d)$ be countable. Set $\Phi=\pi_X^{-1}(\Psi_\d)$. Then $\Lambda_X\rest M_X=\Phi$. 
\end{lemma}
\begin{proof} Let $\U\in M_X$ be such that $\U\in dom(\Phi)\cap dom(\Lambda_X)$. Set $b=\Lambda_X(\U)$. It follows from (2) of Definition \ref{ub strategy} that $b\in M_X$. Because $M_X\models {\sf{gUBH}}$, it follows that $\Phi(\U)=b$. 
\end{proof}
  
\begin{theorem}\label{strategies can be extended}
Suppose $\d$ is a Woodin cardinal and $\eta\geq \d$ is an inaccessible cardinal. Let $g\subseteq Coll(\omega, \eta)$ be generic. Then, in $V[g]$, there is an $Ord$-strategy $\Sigma$\footnote{Recall that we are assuming self-iterability.} for $\V|\d$ such that the following hold.
\begin{enumerate}
\item $\Psi_\d \subseteq \Sigma$,
\item Letting $\Delta$ be the $\omega_1$-fragment of $\Sigma$, $V[g]\models ``\Delta$ is universally Baire". 
\item For all $V[g]$-generic $h$, letting $\Delta^h$ be the canonical extension of $\Delta$ to $V[g*h]$, $\Delta^h\rest V[g]\subseteq \Sigma$. 
\end{enumerate}
 \end{theorem}
 \begin{proof} Let $\l > \eta$ be an inaccessible cardinal. Set $W=W_\l$, $\P=\V|\d$ and given a iteration $\T$ of $\P$ of limit length and a cofinal well-founded branch $b$ of $\T$, set $\psi[\T, b]=  \phi[\T^\frown \{b\}] \wedge  \forall \a<lh(\T) \phi[\T\rest \a+1]$.

Working in $V_\l[g]$, let $\Sigma$ be the strategy given by $\psi$. More precisely, let $\Sigma$ be defined as follows.
\begin{enumerate}
\item $\T\in dom(\Sigma)$ if and only if $lh(\T)$ is of limit length and for every limit $\a<lh(\T)$ if $b=[0, \a)_\T$ then $V_\l[g]\models \psi[\T, b]$.
\item $\Sigma(\T)=b$ if and only if $V_\l[g]\models \psi[\T, b]$. 
\end{enumerate}

The following is an immediate consequence of our definitions. 

\begin{lemma}\label{capturing1} Suppose $X\prec (W |\Psi_\d)$ is countable. Let $k\in V$ be $M_X$-generic. Suppose $(\U, b)\in M_X[k]$ is such that $M_{X}[k]\models \psi[ \U, b]$. Then $\U\in dom(\Lambda_X)$ and $\Lambda_X(\U)=b$.
\end{lemma}

We now work towards showing that $\Sigma$ is a total strategy.
\begin{lemma}\label{step1} Suppose $\T\in dom(\Sigma)$. Then there is at most one branch $b$ such that $V[g]\models \psi[\T, b]$. 
\end{lemma}
 \begin{proof} Towards a contradiction assume not.  Let $X\prec (W|\Psi_\d)$ be countable and $k\subseteq Coll(\omega, \eta_X)$ be $M_X$-generic with $k\in V$. Fix now $\U, b, c\in M_{X}[k]$ such that $M_{X}[k]\models \psi[\U, b]\wedge \psi[\U, c]$. It follows from \rlem{capturing1} that $b=\Lambda_X(\U)=c$. Therefore, $b=c$. 
\end{proof}

\begin{lemma}\label{step2} Suppose $\T\in dom(\Sigma)$. Then there is a branch $b$ such that $V_\l[g]\models \psi[\T, b]$. 
\end{lemma}
  \begin{proof} Towards a contradiction assume not. Let $X\prec (W|\Psi_\d)$ be countable and  $k\subseteq Coll(\omega, \eta_X)$ be $M_X$-generic. It follows that there is an iteration $\U\in M_X[k]$ of $\P_X$ such that\\\\
  (a) for every $\a<lh(\U)$, letting $b_\a=[0, \alpha)_\U$, $M_X[k]\models \psi[\U\rest \a, b_\a]$ but\\
  (b) for no well-founded cofinal branch $b\in M_X[k]$ of $\U$, $M_X[k]\models \psi[\U, b]$. \\\\
  It follows from (a) and \rlem{capturing1} that $\U\in dom(\Lambda_X)$. Hence, setting $\Lambda_X(\U)=b$, $b\in M_X[k]$ and $M_{X}[k]\models \phi[\U^\frown \{b\}]$. Therefore, $M_X[k]\models \psi[\U, b]$.
\end{proof}

\begin{lemma}\label{lemma1 towards uB}  Let $X\prec (W| \Psi_\d)$ be countable and let $k\in V$ be $M_X$-generic for $Coll(\omega, \eta_X)$. Let $\Phi$ be the strategy of $\P_X$ defined by $\psi$ in $M_X[k]$.  Then $\Lambda_X\rest M_X[k]=\Phi$. 
\end{lemma}
\begin{proof}
Suppose that $\T\in M_X[k]$ is according to both $\Lambda_X$ and $\Phi$. Set $b=\Phi(\T)$. Because $\Phi(\T)=b$ we have that $M_X[k]\models \phi[\T^\frown \{b\}]$. Hence, $\Lambda_X(\T)=b$. 
\end{proof}

\begin{corollary}\label{total strategy}
$V_\l[g]\models ``\Sigma$ is a total strategy extending $\Psi_{\d}\rest V_\l$".
\end{corollary}
\begin{proof}
\rlem{step1} and \rlem{step2} imply that $\Sigma$ is a total strategy. To show that it extends $\Psi_\d\rest V_\l$, we reflect. Let $X\prec (W|\Psi_\d)$ be countable and let $k\subseteq Coll(\omega, \eta_X)$ be $M_X$-generic such that $k\in V$. Let $\Phi$ be the strategy of $\P_X$ defined by $\psi$ over $M_X[k]$. It follows from \rlem{lemma1 towards uB} that $\Phi=\Lambda_X\rest (M_X[k])$. It follows from \rlem{simple capturing} that $\Lambda_X\rest M_X=\pi_X^{-1}(\Psi_\d)$.  Hence, $\pi^{-1}_X(\Psi_\d)\subseteq \Phi$.
\end{proof}
 
We now work towards showing that $\Delta=_{def}\Sigma\rest HC^{V[g]}$ is universally Baire. For this it is enough to show that $\psi$ is generically correct. More precisely, it is enough to show that in $V[g]$, for a club of $X\prec (W, \Psi_\d)$ such that $V_\eta\cup \{\eta\}\subseteq X$, whenever $k\in V[g]$ is $M_X[g]$-generic and $(\T, b)\in M_X[g][k]$,
\begin{center}
$M_X[g][k]\models \psi[\T, b] \iff V[g]\models \psi[\T, b]$.\footnote{See \cite[Lemma 4.1]{DMT} for a proof of the equivalence.}
\end{center}

Working in $V$, fix $X\prec H_{\l^+}$ such that $W, \Psi_\d\in X$. It is enough to show that our claim holds in $M_X$. Let $k\in V$, $k\subset Coll(\omega, \eta_X)$ be $M_X$-generic. Let $\Phi$ be the strategy defined by $\psi$ over $M_X[k]$ and $\Psi=\pi_X^{-1}(\Psi_\d)$. Let $Y\prec (W_X| \Psi)$ be any countable substructure in $M_X[k]$ such that $V^{M_X}_{\eta_X}\cup \eta_X \subset Y$ and let $h\in M_X[k]$ be $M_Y[k]$-generic. Fix $(\T, b)\in M_Y[k][h]$.

Suppose now that $M_Y[k][h]\models \psi[\T, b]$. Because $\pi_X[Y]\in V$ we have that $\T$ is according to $\Lambda_Y$ and $\Lambda_Y(\T)=b$. But because $\pi_Y\rest \eta_X=id$, we have that $\Lambda_Y=\Lambda_X$. Therefore, $\T$ is according to $\Lambda_X$ and $\Lambda_X(\T)=b$. It follows from \rlem{lemma1 towards uB} that $\Phi(\T)=b$, i.e. $M_X[k] \vDash \psi[\T,b]$. The reader can easily verify that these implications are reversible, and so if $\Phi(\T)=b$ then $M_Y[k][h]\models \psi[\T, b]$.

Finally, we need to verify that if $h$ is $V[g]$-generic for a poset of size $<\l$ then $\Delta^h\rest V_\l[g]\subseteq \Sigma$. This again can be verified by first reflecting in $V$. Indeed, working in $V$, fix $X\prec H_{\l^+}$ be countable such that $W, \Psi_\d\in X$. Let $(k, \Phi, \Psi)$ be as above. Let $\Gamma=\Phi\rest HC^{M_X[k]}$. Let $h\in V$ be any $M_X[k]$-generic. We want to see that $\Gamma^h\rest M_X[k]\subseteq \Phi$. To see this, let $\T\in M_X[k]$ be according to both $\Gamma^h$ and $\Phi$. Let $b=\Gamma^h(\T)$. It follows that $M_X[k][h]\models \psi[\T, b]$. Hence, $\T\in dom(\Lambda_X)$ and $\Lambda_X(\T)=b$. It follows from \rlem{lemma1 towards uB} that $\Phi(\T)=b$. 

Thus far we have shown that \rthm{strategies can be extended} holds in $V_\l[g]$ for any inaccessible $\l>\eta$. Let $\Sigma_\l$ be the strategy defined above. To finish the proof of \rthm{strategies can be extended} it is enough to show that if $\l_0<\l_1$ are two inaccessible cardinals bigger than $\eta$ then $\Sigma_{\l_1}\rest V_{\l_0}[g]=\Sigma_{\l_0}$. This can be verified by a reflection argument similar to the ones given above.

Indeed, let $X\prec H_{\l_1^+}$ be countable such that $\{W_{\l_0}, W_{\l_1}\} \in X$. Let $k\subseteq Coll(\omega, \eta_X)$ be $M_X$-generic such that $k\in V$. Let $\Phi_0$ and $\Phi_1$ be the versions of $\Sigma_{\l_0}$ and $\Sigma_{\l_1}$ in $M_X[k]$. It follows from \rlem{lemma1 towards uB} that for $i\in 2$, $\Phi_i=\Lambda_X\rest M_{X\cap W_{\l_i}}[k]$. Therefore, $\Phi_0\subseteq \Phi_1$. This completes the proof of Theorem \ref{strategies can be extended}.
\end{proof}

We record a useful corollary to the proof of \rthm{strategies can be extended}. We let $\psi$ be the formula used in the proof of \rthm{strategies can be extended}. If $g, \Sigma$ are as in \rthm{strategies can be extended} and $k$ is $V[g]$-generic then we let $\Sigma^k$ be the extension of $\Sigma$ to $V[g][k]$. 

\begin{corollary}\label{definable extension} Suppose $\delta, g, \Sigma$ are as in Theorem \ref{strategies can be extended}. Suppose $\l$ is an inaccessible cardinal and $k$ is $V[g]$-generic for a poset in $V_\l[g]$. Then $\Sigma^k\rest V_\l[g][k]$ is defined via $\psi$. More precisely, the following conditions hold.
\begin{enumerate}
\item $\T\in dom(\Sigma^k)\cap V_\l[g*k]$ if and only if for every limit $\a<lh(\T)$, setting $b_\a=[0, \a)_\T$, $V_\l[g*k]\models \psi[\T\rest \a, b_\a]$.
\item For $\T\in dom(\Sigma^k)\cap V_\l[g*k]$, $\Sigma^k(\T)=b$ if and only if $V_\l[g*k]\models \psi[\T, b]$. 
\end{enumerate}
\end{corollary}

As the definition of $\Sigma$ uses only parameters from $V$, it follows that in all generic extensions $V[h]$ of $V$, $\Psi_\d$ has an extension $\Psi_\d^h$.  For instance, we can define $\Psi^h_\d(\U)$ by first selecting some inaccessible $\eta$ such that $h$ is generic for a poset in $V_\eta$ and $\U\in V_\eta[h]$ then picking a generic $g\subseteq Coll(\omega, \eta)$ such that $V[h]\subseteq V[g]$ and then finally setting $\Psi^h_\d(\U)=\Sigma(\U)$ where $\Sigma$ is as in \rthm{strategies can be extended}. 

\section{Some correctness results}

Say $u=(\eta,\d, \l)$ is a \textit{good triple} if it is increasing, $\d$ is a Woodin cardinal, and $\l$ is an inaccessible cardinal. The assumption on Woodinness of $\d$ will not be necessary in this section but will be used extensively in subsequent sections. Fix a good triple and set $\Phi=\Psi_\d\rest H_\l$. The goal of this section is to show that many Skolem hulls of $\Phi$ are computed correctly. Forcing posets in some of the main claims this section will be in $V_\eta$.  We start by showing that a stronger form of \rlem{simple capturing} holds.

\begin{lemma}\label{first step}  Suppose $X\prec ((W_\l, u) | \Phi)$ is countable and $k\in V$ is $M_X$-generic. Then 
\begin{center}
$\Phi_X^k\rest (M_X[k])=\Lambda_X\rest (M_X[k])$\footnote{Here $\Phi^k_X$ is the generic interpretation of $\Phi_X$ in $M_X[k]$ using the definition of $\Phi$ given in \rthm{strategies can be extended}.}.
\end{center} 
\end{lemma}
\begin{proof} Fix $\T\in dom(\Phi_X^k)\cap dom(\Lambda_X)$ and set $\Phi_X^k(\T)=b$. It follows from \rcor{definable extension} that $\M_X[k]\models \psi[\T, b]$. Therefore, $\Lambda_X(\T)=b$.  
\end{proof}

The following is a straightforward corollary of \rlem{first step} and can be proven by a reflection like that in the proof of Theorem \ref{strategies can be extended}. 

\begin{corollary}\label{first step generically} Suppose $g$ is generic for a poset in $V_\eta$ and $X\prec ((W_\l, u) | \Phi^g)$ is countable in $V[g]$. Let $k\in V[g]$ be $M_X$-generic. Then 
\begin{center}
$\Phi_X^k\rest (M_X[k])=\Lambda_X\rest (M_X[k])$.
\end{center} 
\end{corollary}

\begin{corollary}\label{generic interpretability claim} Suppose $g$ is generic for a poset in $V_\eta$ and $i: \V\rightarrow \P$ is an iteration embedding via a normal iteration $\T$ of length $<\l$ that is based on $\V|\d$ and is according to $\Phi$. Then $i(\Phi)=\Phi^g_{\P|i(\d)}\rest \P$. \footnote{Recall that $\Phi_{\P|i(\d)} =_{\textrm{def}} \Phi_{\T, \P|i(\d)}$ is the tail strategy of $\P|i(\delta)$ induced by $\Phi$.}
\end{corollary}
\begin{proof} It is enough to prove the claim in some $M_Z$ where $Z\prec ((H_{\l^{++}}, W_\l, u, \Phi)|\Phi)$ is countable. Let $h\in V$ be $M_Z$-generic for a poset in $M_Z|\eta_Z$, and let $\U\in M_Z|\l_Z[h]$ be a normal iteration of $\V_Z$ based on $\V_Z|\delta_Z$ according to $\Phi_Z$ with last model $\Q$. We want to see that $\pi^\U(\Phi_Z)=(\Phi_Z)_{\Q|\pi^\U(\d_Z)}^h\rest \Q$. 

Let $\R$ be the last model of $\pi_Z\U$ and $\sigma:\Q\rightarrow \R$ come from the copying construction. It follows from \cite[Theorem 4.9.1]{Neeman} that $\sigma$ is generic over $\R$ and $\R[\sigma]\in V$. It then follows from \rcor{first step generically} that $\pi^\U(\Phi_Z)=(\pi^{\pi_Z\U}(\Phi))^\sigma$. It again follows from \rcor{first step generically}  that $\Phi_Z^h=\Lambda_Z\rest M_Z[h]$, and hence 

\begin{center}
$(\Phi_Z)_{\Q|\pi^\U(\d_Z)}^h\rest \Q=(\Lambda_Z)_{\Q|\pi^\U(\d_Z)}\rest \Q = (\pi^{\pi_Z\U}(\Phi))^\sigma = \pi^\U(\Phi_Z)$.
\end{center}
\end{proof}

\begin{figure}
\centering
\begin{tikzpicture}[node distance=2.5cm, auto]
  \node (A) {$\V$};
  \node (B) [right of=A] {$\P$};
  \node (C) [below of=A] {$M_Z$};
  \node (D) [right of=C] {$\Q$};
   \draw[->] (A) to node {$i$}(B); 
  \draw[->] (C) to node {$\sigma$}(D);
  \draw[->] (C) to node  {$\pi_Z$} (A);
  \draw[->] (D) to node {$\tau$} (B);
  \end{tikzpicture}
\caption{Corollary \ref{generic interpretability}}
\label{fig:comm_diag}
\end{figure}

\begin{corollary}[Figure \ref{fig:comm_diag}]\label{generic interpretability} Suppose $g$ is generic for a poset in $V_\eta$ and $i: \V\rightarrow \P$ is an iteration embedding via a normal iteration $\T$ of length $<\l$ that is based on $\V|\d$ and is according to $\Phi$. Let $X\prec ((W_\l, u) | \Phi^g)$ be countable in $V[g]$ and let $\Q\in HC^{V[g]}$ be such that there are embeddings $\sigma: M_X\rightarrow \Q$  and $\tau:\Q\rightarrow\P$ with the property that  $i\circ \pi_X=\tau\circ \sigma$. Then for any $\Q$-generic $k\in V[g]$,
\begin{center}
$(\sigma(\Phi_X))^k_{\Q|\sigma(\delta_Z)}=(\tau$-pullback of $\Phi^g_{\P|i(\d)})\rest \Q[k]$.
\end{center}
\end{corollary}
\begin{proof} It is enough to prove the claim assuming $g$ is trivial. The more general claim then will follow by using the proof of \rcor{generic interpretability claim}. It follows from \cite[Corollary 4.9.2]{Neeman} that $\tau$ is generic over $\P$ and that $\P[\tau]$ is a definable class of $V$; here, to apply \cite[Corollary 4.9.2]{Neeman}, we need $\Q$ is countable in $V[g]$. Applying \rcor{first step generically} and \rcor{generic interpretability claim} in $\P$, we get that 
\begin{center}
$(\sigma(\Phi_X))^k_{\Q|\sigma(\delta_Z)}=(\tau$-pullback of $\Phi^g_{\P|i(\d)})\rest \Q[k]$.
\end{center}
\end{proof}

\begin{corollary}\label{generic interpretability claim1} Suppose $i: \V\rightarrow \P$ is an iteration embedding via a normal iteration $\T$ of length $<\l$ that is based on $\V|\d$ and is according to $\Phi$. Let $h\in V$ be $\P$-generic for a poset in $\P|\l$. Then $i(\Phi)^h=\Phi_{\P|i(\d)}\rest \P[h]$. 
\end{corollary}
\begin{proof} It is enough to prove the claim in some $M_Z$ where $Z\prec ((H_{\l^{++}}, W_\l, u, \Phi)|\Phi)$. Let $\U\in M_Z$ be an iteration of length $< \lambda_Z$ on $M_Z$ based on $M_Z|\delta_Z$ and $j: M_Z\rightarrow\Q$ be the iteration embedding. Let $G\in M_Z$ be $\Q$-generic for a poset in $\Q|\l_Z$.  We want to see that $j(\Phi_Z)^{G}=(\Phi_Z)_{\Q|j(\d_Z)}\rest \Q[G]$.  Let $\T=\pi_Z\U$, $\P$ the last model of $\T$, $k=\pi^\T$ and $\tau:\Q\rightarrow \P$ be the copy map.

It follows from \rcor{generic interpretability} that $j(\Phi_Z)^G=(\tau$-pullback of $\Phi_{\P|i(\d)})\rest \Q[G]$. But because $\Phi_Z=\Lambda_Z\rest M_Z$ (see \rlem{first step}), 
\begin{center}
$(\tau$-pullback of $\Phi_{\P|i(\d)})\rest \Q[G]=(\Phi_Z)_{\Q|j(\d_Z)}\rest \Q[G]$.
\end{center}
Therefore,
\begin{center}
$j(\Phi_Z)^G = (\Phi_Z)_{\Q|j(\delta_Z)}\rest \Q[G]$.
\end{center}
\end{proof}

Suppose $M$ is a transitive model of set theory and $\nu$ is its least strong cardinal.  Suppose $M\models ``u=(\eta,\d, \l)$ is a good triple" and suppose $\T$ is a normal iteration of $M$. We say $\T$ is a \textit{sealed} iteration if $\T=\T_0^\frown\{E_0\}$ is such that 
\begin{enumerate}
\item $\T_0$ is a normal iteration of $M$ of successor length based on $M|{\d}$ with last model $N$, 
\item $\T_0$ is above $\nu$ (this implies that $\delta > \nu$),
\item $E_0\in N$ is an extender such that $\cp(E_0)=\nu$, $lh(E_0)>\pi^{\T_0}(\d)$,
\item $N$ has an inaccessible cardinal in the interval $(\pi^{\T_0}(\d), lh(E_0))$.
\end{enumerate}
Clearly the last model of $\T$ is $Ult(M, E_0)$. We say that a normal iteration $\T$ is a stack of sealed iterations if for some $n<\omega$, $\T=\oplus_{i\leq n}\T_i$ such that $\T_i$ is a sealed iteration of its first model. 

\begin{corollary}\label{coherence} Suppose $u=(\eta, \d, \l)$ is a good triple, $g$ is generic for a poset in $V_\eta$ and $\T\in V_\l[g]$ is a normal iteration of $\V$ that is a stack of sealed iterations and is according to $\Phi^g$ where $\Phi=\Psi_\d$. Set $\T=\oplus_{i\leq n}\T_i$ and let $\P$ be the last model of $\T_{n-1}$ if $n>0$ and $\V$ otherwise. Let $\T_n=(\U, E)$ and let $\Q$ be the last model of $\U$. Set $\nu=\pi^{\U}(\pi^{\oplus_{i<n}\T_i}(\d))$. Then $\Phi^g_{Ult(\P, E)|\nu}=\Phi^g_{\Q|\nu}$. 
\end{corollary}
\begin{proof} We prove the claim in some $M_Z$ where $Z\prec ((H_{\l^+}, W_\l, u, \Phi)|\Phi)$ is countable. Let $h$ be $M_Z$-generic for a poset in $M_Z|\eta_Z$ and let $(\W, \R, \W_n, \S, \X,  F, \xi)\in M_Z[h]$ play the role of $(\T, \P, \T_n, \Q, \U, E, \nu)$. 

We will redefine the objects $\P$ etc. in the following; this will not cause any confusion as we have no more use for the original objects. Let $\P$ be the last model of the $\pi_Z$-copy of $\oplus_{i<n}\W_i$ and let $\sigma:\R\rightarrow \P$ be the copy map. We have that $\sigma$ is generic over $\P$ (see \cite[Corollary 4.9.2]{Neeman}) and $\P[\sigma]$ is a definable class of $V$. Let $\Q$ be the last model of $\sigma\X$ and let  $\tau_0:\S\rightarrow \Q$ and $\tau_1:Ult(\R, F)\rightarrow Ult(\P, \tau_0(F))$ come from the copying construction. Notice that 
\begin{center}
$\tau_0\rest (\S|lh(F))=\tau_1\rest (\S|lh(F))$.
\end{center}
We then let $\tau$ be this common embedding. Set $\tau_0(F)=E$ and $\nu=\tau_0(\xi)$. We have that $\tau_0$ and $\tau_1$ are generic over $\Q$ and $Ult(\P, E)$ respectively. 

We now want to see that in $M_Z[h]$, 
\begin{center}
$(\Phi_Z^h)_{Ult(\R, F)|\xi}=(\Phi_Z^h)_{\S|\xi}$. 
\end{center}
Notice that it follows from \rlem{first step} that $\Phi^h_Z=\Lambda_Z\rest M_Z[h]$. Let $\Gamma_0=(\tau$-pullback of $\Phi_{\Q|\nu})$ and $\Gamma_1=(\tau$-pullback of $\Phi_{Ult(\P, E)|\nu})$. It follows that \\\\
(0) $\Gamma_0\rest M_Z[h]=(\Phi_Z^h)_{\S|\xi}$ and $\Gamma_1\rest M_Z[h]=(\Phi_Z^h)_{Ult(\R, F)|\xi}$.\\\\
Let $i:\V\rightarrow \Q$ and $j: \V\rightarrow Ult(\P, E)$ be the iteration maps. It follows from \rcor{generic interpretability claim} that\\\\
(1) $\Phi_{\Q|\nu}\rest \Q=i(\Phi)_{\Q|\nu}$ and $\Phi_{Ult(\P, E)|\nu}=j(\Phi)_{Ult(\P,E)|\nu}$.\\\\
Because $Ult(\P, E)|lh(E)=\Q|lh(E)$, $lh(E)>\nu$ is an inaccessible cardinal in $\Q$, and $\Q|lh(E)\models \sf{gUBH}$, we have that\\\\
(2) $i(\Phi)_{\Q|\nu}\rest (\Q|lh(E))=j(\Phi)_{Q|\nu}\rest (\Q|lh(E))=_{def}\Sigma$\\\\
implying by the way of (1) that\\\\
(3) $\Phi_{\Q|\nu}\rest (\Q|lh(E))=\Phi_{Ult(\P, E)|\nu}\rest (\Q|lh(E))$.\\\\
 Using   \cite[Corollary 4.9.2]{Neeman} we can find $H\in V$ that is $\Q$-generic for a poset in $\Q|\nu$ and is such that $\tau_0\in \Q[H]$. It now follows that $\tau\in Ult(\P, E)[H]$ as $\tau\in \Q|lh(E)[H]$. We now have  that\\\\
 (4) $(\Sigma^H)^{Ult(\P, E)[H]}\rest (\Q|lh(E)[H])=(\Sigma^H)^{\Q[H]}\rest (\Q|lh(E)[H])$.\\\\
 Applying \rcor{generic interpretability claim1} to (4) we get that\\\\
 (5) $\Phi_{\Q|\nu}\rest (\Q|lh(E)[H])=\Phi_{Ult(\P, E)|\nu}\rest (\Q|lh(E)[H])$.\\\\
It follows from (5) that \\\\
(6) $\Gamma_0\rest M_Z[h]$ and $\Gamma_1\rest M_Z[h]$ are equal.\\\\
(6) then implies, by the way of (0), that $(\Phi_Z^h)_{Ult(\R, F)|\xi}=(\Phi_Z^h)_{\S|\xi}$.
\end{proof}

\section{Capturing universally Baire sets}\label{sec: cap ub}

The following is a useful corollary of \rthm{strategies can be extended}. We say that a pair of trees $T,S$ are \textit{$\delta$-absolutely complementing} if for any poset $\mathbb{P}$ of size $\leq \d$, for any generic $g\subseteq \mathbb{P}$, $V[g]\models ``p[T]=\bR-p[S]"$. Similarly, we say that $T,S$ are \textit{$<\delta$-absolutely complementing} if for any poset $\mathbb{P}$ of size $< \d$, for any generic $g\subseteq \mathbb{P}$, $V[g]\models ``p[T]=\bR-p[S]"$. Given a limit of Woodin cardinals $\nu$ and $g\subseteq Coll(\omega, <\nu)$, let
\begin{enumerate}
\item $\bR^*_g=\bigcup_{\a<\nu}\bR^{V[g\cap Coll(\omega, \a)]}$, 
\item $\Delta_g$ be the set of reals $A\in V(\bR^*)$ such that for some $\a<\nu$, there is a pair $(T, S)\in V[g\cap Coll(\omega, \a)]$ such that $V[g\cap Coll(\omega, \a)]\models ``(T, S)$ are $<\nu$-complementing trees" and $p[T]^{V(\bR^*)}=A$, and
\item $DM(g)=L(\Delta_g, \bR^*_g)$. 
\end{enumerate}
 
 The following is immediate from results of the previous sections. 
 \begin{corollary}\label{strategy in der model} Suppose $\nu$ is a limit of Woodin cardinals. Let $\d<\nu$ be a Woodin cardinal, and let $g\subseteq Coll(\omega, <\nu)$ be $V$-generic. Then $\Psi^g_\d\in DM(g)$. 
 \end{corollary}
 
 We next need a characterization of universally Baire sets via strategies. We show this in \rlem{ub to strategy}. The lemma is standard. 
 
  If $\nu$ is a Woodin cardinal we let $\sf{EA}_\nu$ be the $\omega$-generator version of the extender algebra associated with $\nu$ (see e.g. \cite{steel2010outline} for a detailed discussion of Woodin's extender algebras).  We say the triple $(M, \d, \Phi)$ \textit{Suslin, co-Suslin captures}\footnote{This notion is probably due to Steel, see \cite{DMATM}.} the set of reals $B$ if there is a pair $(T, S)\in M$ such that $M\models ``(T, S)$ are $\d$-complementing" and 
 \begin{enumerate}
 \item $M$ is a countable transitive model of some fragment of $\sf{ZFC}$,
 \item $\Phi$ is an $\omega_1$-strategy for $M$,
 \item $M\models ``\d$ is a Woodin cardinal",
 \item for $x\in \bR$, $x\in B$ if and only if there is an iteration $\T$ of $M$ according to $\Phi$ with last model $N$ such that $x$ is generic over $N$ for  $\textsf{EA}_{\pi^\T(\d)}^{N}$ and $x\in p[\pi^\T(T)]$.
 \end{enumerate}
 
 The next lemma is standard and originates in \cite{IT}.

  \begin{lemma}\label{realizability} Suppose $u=(\eta, \d, \l)$ is a good triple and $g$ is $V$-generic for a poset in $V_\eta$. Suppose $X\prec (W_\l[g]|\Psi_{\eta, \d}^g)$ is countable in $V[g]$. Then whenever $\T$ is a countable iteration of $M_X$ according to $\Lambda_X$ with last model $N$, there is $\sigma:N\rightarrow W_\l[g]$ such that $\pi_X=\sigma\circ \pi^\T$.  
 \end{lemma}
 \begin{proof} Let $\P = W_\lambda[g]$. Let $\U =_{def} \pi_X\T$ be the copy of $\T$, considered as a tree on $V[g]$. Let $W$ be the last model of $\U$.  There is then $\tau:N\rightarrow \pi^{\U}(\P)$ such that $\pi^\U\circ \pi_X=\tau\circ \pi^\T$. It follows by absoluteness, noting $N\in W$ is countable and $\pi^\U(\P)\in W$, that there is $m:N\rightarrow \pi^\U(\P)$ with $m\in W$ such that $\pi^\U(\pi_X)=m\circ \pi^\T$. The existence of $\sigma$ follows from elementarity. 
 \end{proof}
 
 The next lemma is also standard, but we do not know its origin. To state it we need to introduce some notations. Suppose $M$ is a countable transitive model of set theory and $\Phi$ is a strategy of $M$. Let $(\eta, g)$ be such that $g$ is $M$-generic for a poset in $M|\eta$. Let $\Phi'$ be the fragment of $\Phi$ that acts on iterations that are above $\eta$. Then $\Phi'$ can be viewed as an iteration strategy of $M[g]$. This is because if $\T$ is an iteration of $M[g]$ above $\eta$, there is an iteration $\U$ of $M$ that is above $\eta$ and such that 
\begin{enumerate}
\item $lh(\T)=lh(\U)$,
\item $\T$ and $\U$ have the same tree structure,
\item for each $\a<lh(\T)$, $M^\T_\a=M^\U_\a[g]$,
\item for each $\a<lh(\T)$, $E_\a^\T$ is the extension of $E_\a^\U$ onto $M_\a^\U[g]$.
\end{enumerate}
Let $\Phi''$ be the strategy of $M[g]$ with the above properties. We then say that $\Phi''$ is induced by $\Phi'$. We will often confuse $\Phi''$ with $\Phi'$.

\begin{corollary}\label{strategy for the generic extension} Suppose $(\eta, \d, \l)$ is a good triple, $g$ is generic for a poset of size $<\eta$ and $h\subseteq Coll(\omega, \l)$ is generic over $V$ such that $V[g]\subseteq V[h]$. Let $\Sigma$ be as in \rthm{strategies can be extended} applied to $h$ and $\Psi_\d$, and let $\Phi$ be the fragment of $\Sigma\rest V[g]$ that acts on iterations that are above $\eta$. Then $\Phi$ induces a strategy $\Phi'$ for $\V|\d[g]$, and $\Phi'$ is projective in $\Phi$. \footnote{This just means in $V[h]$, $\Phi'\rest HC$ is definable over the structure $(HC, \in, \Phi\rest HC)$ perhaps with parameters in $HC$.}
\end{corollary}

We can now state our lemma.
 
 \begin{lemma}\label{ub to strategy} Suppose $u=(\eta, \d, \l)$ is a good triple and $g$ is $V$-generic for a poset in $V_\eta$. Let $A\in \Gamma^\infty_g$.  Then, in $V[g]$, there is a club of countable $X\prec (W_\l[g]| \Psi_{\eta, \d}^g)$ such that $(M_X, \d_X, \Lambda^g_X)$ Suslin, co-Suslin captures $A$.\footnote{To conform with the above setup, we tacitly assume $\Lambda^g_X$ to be the iteration strategy acting on trees above $\eta_X$.} For each such $X$, let $X' = X\cap W_\lambda \prec W_\lambda$, and $(M_{X'}, \Lambda_{X'})$ be the transitive collapse of $X'$ and its strategy. Then $A$ is projective in $\Lambda_{X'}$. Moreover, these facts remain true in any further generic extension by a poset in $V_\eta[g]$. 
 \end{lemma}
 \begin{proof} 
Let $\P = W_\lambda[g]$. Work in $V[g]$. Let $(T, S)$ be $\l$-complementing trees such that $A=p[T]$. Let $X\prec W_\l[g]$ be countable such that $(T, S)\in X$. We claim that $(M_X, \d_X, \Lambda^g_X)$ Suslin co-Suslin captures $A$. Let $\delta = \delta_X$. To see this fix a real $x$. Let $\T$ be any countable normal iteration of $M_X$ such that
 \begin{enumerate}
 \item $\T$ is according to $\Lambda^g_X$,
 \item $\T$ has a last model $N$,
 \item $x$ is generic for ${\sf{EA}}_{\pi^\T(\d)}^N$.
 \end{enumerate}
 Using \rlem{realizability}, we can find $\sigma:N\rightarrow \P$ such that $\pi_X=\sigma\circ \pi^\T$. 
 
 Assume first $x\in A$. Then $x\in p[T]$. If now $x\not \in p[\pi^\T(T_X)]$ then $x\in p[\pi^\T(S_X)]$ (this uses the fact that $T_X,S_X$ are $\lambda_X$-complementing in $M_X$) and hence, $x\in p[S]$ (this follows from the fact that $\sigma[\pi^\T(S_X)]\subseteq S$). Thus, $x\in p[\pi^\T(T_X)]$. 
 
 Next suppose $x \in p[\pi^\T(T_X)]$. Then because $\sigma[\pi^\T(T_X)]\subseteq T$, $x\in p[T]$ implying that $x\in A$. 
 
 That $\Lambda^g_X$ is projective in $\Lambda_{X'}$ follows from \rcor{strategy for the generic extension}; hence $A$ is projective in $\Lambda_{X'}$. We leave it to the reader to verify that these facts remain true in a further generic extension by a poset in $V_\eta[g]$.
 \end{proof}

 \section{A derived model representation of $\Gamma^\infty$}\label{sec:der model}
 
 In this section our goal is to establish a derived model representation of $\Gamma^\infty$. We set $\iota=\kappa^+$ and fix $g\subseteq Coll(\omega, \iota)$.

We say $u=(\eta, \d, \d', \l)$ is a \textit{good quadruple} if $(\eta, \d, \l)$ and $(\eta, \d', \l)$ are good triples with $\d<\d'$.  Suppose $u=(\eta, \d, \d', \l)$ is a good quadruple and $h$ is a $V[g]$-generic such that $g*h$ is generic for a poset in $V_\eta$.  Working in $V[g*h]$, let $D(h, \eta, \d, \l)$ be the club of countable 
\begin{center}
$X\prec ((W_\l[g*h], u)| \Psi^g_{\eta, \d})$
\end{center}
 such that $H_{\iota}^V\cup\{g\}\subseteq X$. 

 Suppose $A\in \Gamma^\infty_{g*h}$. Then for a club of $X\in D(h, \eta, \d, \l)$, $A$ is Suslin, co-Suslin captured by $(M_X, \d_X, \Lambda^{g*h}_X)$ and $A$ is projective in $\Lambda_{X'}$ where $X'  =X\cap W_\l$ (see \rlem{ub to strategy}). Given such an $X$, we say $X$ \textit{captures} $A$.

Let $k\subseteq Coll(\omega, \Gamma^\infty_{g*h})$ be generic, and let $(A_i:i<\omega),(w_i: i<\omega)$ be generic enumerations of $\Gamma^\infty_{g*h}$ and $\bR_{g*h}$ respectively in $V[g*h*k]$. Let $(X_i: i<\omega)\in V[g*h*k]$ be such that for each $i$
\begin{enumerate}
\item $X_i \in D(h,\eta,\delta,\lambda)$, and
\item $X_i$ captures $A_i$. 
\end{enumerate}
 In particular, $A_i$ is projective in $\Lambda_{X^{'}_i}$, where $X^{'}_i = X_i\cap W_\lambda$.
We set  $M^0_n=M_{X'_n}$, 
 $\pi^0_n=\pi_{X_0}$, 
$\k_0=\k_{X_0}$,
 $\nu_0=\d_{X_0}$,
 $\nu_0'=\d_{X_0}'$,
$\eta_0=\eta_{X_0}$,
 $\d_0=\d$,
$\P_0=\V$.

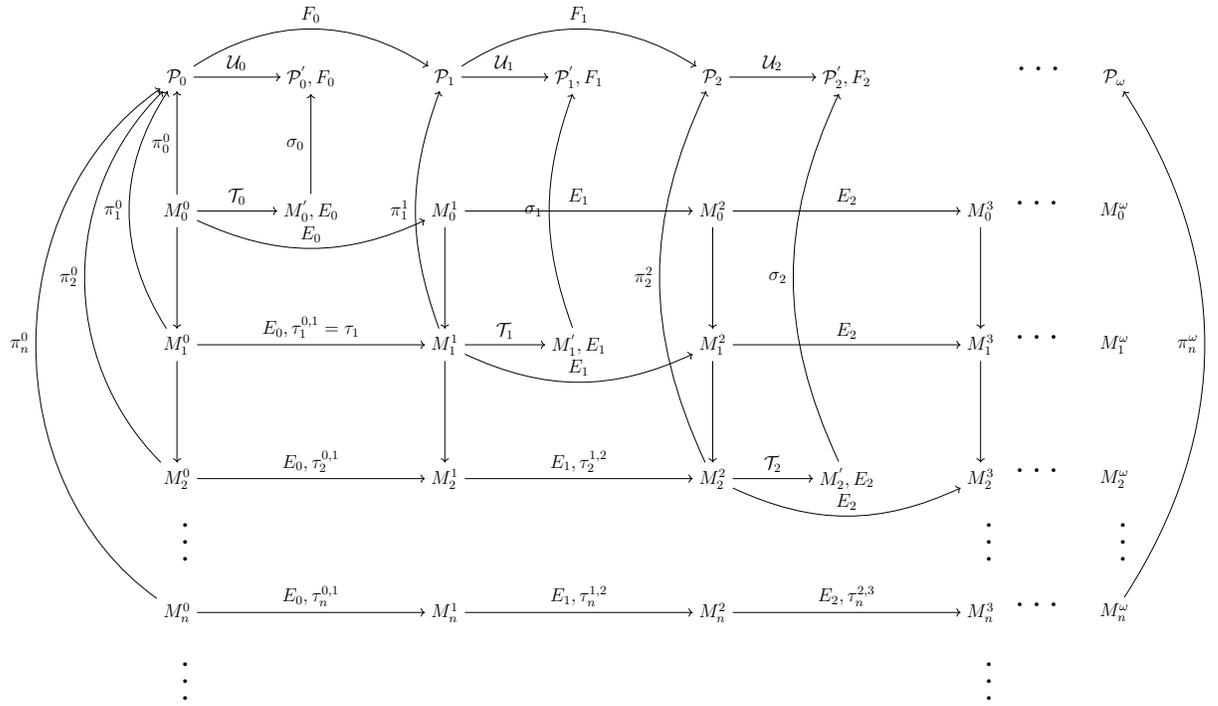
\begin{figure}
\centering
\resizebox{0.85\textheight}{!}{
\begin{tikzpicture}[align=center, node distance= 3cm, auto]
\node (A) {$\P_0$};
\node (B) [right of=A] {$\P^{'}_{0},F_0$};
\draw[->] (A) to node {$\U_0$} (B);
\node (C) [right of=B] {$\P_1$};
\draw[->, bend left=33] (A) to node {$F_0$} (C);
\node (D) [right of = C] {$\P^{'}_{1}, F_1$};
\draw[->] (C) to node  {$\U_1$} (D);
\node (E) [right of=D] {$\P_2$};
\draw [->, bend left=33] (C) to node {$F_1$} (E);
\node (F) [right of=E] {$\P^{'}_{2}, F_2$};
\draw[->] (E) to node {$\U_2$} (F);
\node (G) [right of=F] {};
\node (H) [right of=F]{};
\node (K) [right of=H]{$\P_\omega$};
\path (H) -- (K) node [font=\Huge, midway, sloped] {$\dots$};

\node (AA) [below of=A]{$M^0_0$};
\node (AB) [right of=AA]{$M^{'}_{0}, E_0$};
\node (AC) [right of=AB]{$M^1_0$};
\node (AD) [below of=D]{};
\node (AE) [below of=E]{$M^2_0$};
\node (AF) [right of=AE]{};
\node (AG) [right of=AF]{$M^3_0$};
\node (AH) [right of=AF]{};
\node (AK) [below of=K]{$M^\omega_0$};
\path (AH)--(AK) node [font=\Huge, midway, sloped]{$\dots$};
\draw[->] (AA) to node {$\T_0$} (AB);
\draw[->] (AA) to node {$\pi^0_0$} (A);
\draw[->] (AB) to node {$\sigma_0$}(B);
\draw[->, bend right=25] (AA) to node {$E_0$}(AC);
\draw[->] (AC) to node {$E_1$}(AE);

\node (BA) [below of=AA]{$M^0_1$};
\node (BB) [below of=AB]{};
\node (BC) [right of=BB]{$M_1^1$};
\node (BD) [right of=BC]{$M^{'}_{1}, E_1$};
\node (BE) [right of=BD]{$M^2_1$};
\node (BF) [right of=BE]{};
\node (BG) [below of=AG]{$M^3_1$};
\node (BH) [right of=BF]{};
\node (BK) [right of=BH]{$M^\omega_1$};
\path (BH)--(BK) node [font=\Huge, midway, sloped]{$\dots$};
\draw[->] (BA) to node {$E_0, \tau^{0,1}_1=\tau_1$} (BC);
\draw[->] (BC) to node {$\T_1$} (BD);
\draw[->, bend left=33] (BA) to node {$\pi^0_1$}(A);
\draw[->, bend left=20] (BC) to node {$\pi^1_1$}(C);
\draw[->] (AA) to node {} (BA);
\draw[->] (AC) to node {} (BC);
\draw[->, bend left=20] (BD) to node {$\sigma_1$} (D);
\draw[->, bend right=25] (BC) to node {$E_1$}(BE);
\draw[->] (AE) to node {}(BE);

\node (CA) [below of=BA]{$M^0_2$};
\node (CB) [below of=BB]{};
\node (CC) [below of=BC]{$M^{1}_{2}$};
\node (CD) [right of=CC]{};
\node (CE) [right of=CD]{$M^2_2$};
\node (CF) [right of=CE]{$M^{'}_{2}, E_2$};
\node (CG) [right of=CF]{$M^3_2$};
\node (CH) [right of=CF]{};
\node (CK)[below of=BK]{$M^\omega_2$};
\draw[->] (BA) to node {}(CA);
\draw[->, bend left=45] (CA) to node {$\pi^0_2$}(A);
\draw[->] (CA) to node {$E_0, \tau^{0,1}_2$}(CC);
\draw[->] (CC) to node {$E_1, \tau^{1,2}_2$}(CE);
\draw[->] (CE) to node {$\T_2$}(CF);
\draw[->, bend left=25] (CE) to node {$\pi^2_2$}(E);
\draw[->, bend left=25] (CF) to node {$\sigma_2$}(F);
\draw[->,bend right=25] (CE) to node {$E_2$}(CG);
\draw[->] (BC) to node {}(CC);
\draw[->] (BE) to node {}(CE);
\draw[->] (BE) to node {$E_2$}(BG);
\draw[->] (AE) to node {$E_2$}(AG);
\draw[->] (AG) to node {}(BG);
\draw[->] (BG) to node {}(CG);
\path (CH)--(CK) node [font=\Huge, midway, sloped]{$\dots$};
\node (DA) [below of=CA]{$M^0_n$};
\node (DB) [right of=DA]{};
\node (DC) [right of=DB]{$M^1_n$};
\node (DD) [right of=DC]{};
\node (DE) [right of=DD]{$M^2_n$};
\node (DF) [right of=DE]{};
\node (DG) [right of=DF]{$M^3_n$};
\node (DH) [right of=DF]{};
\node (DK) [right of=DH]{$M^\omega_n$};

\path (CA) -- (DA) node [font=\Huge, midway, sloped] {$\dots$};
\draw[->] (DA) to node {$E_0, \tau^{0,1}_n$}(DC);
\draw[->] (DC) to node {$E_1,\tau^{1,2}_n$}(DE);
\draw[->] (DE) to node {$E_2, \tau^{2,3}_n$}(DG);
\path (DH) -- (DK) node [font=\Huge, midway, sloped] {$\dots$};
\draw[->, bend right=33] (DK) to node {$\pi^\omega_n$}(K);
\node (EA) [below of=DA]{};
\node (EG) [below of=DG]{};
\path (DA) -- (EA) node [font=\Huge, midway, sloped] {$\dots$};
\path (DG) -- (EG) node [font=\Huge, midway, sloped] {$\dots$};
\path (CG) -- (DG) node [font=\Huge, midway, sloped] {$\dots$};
\path (CK) -- (DK) node [font=\Huge, midway, sloped] {$\dots$};
\draw[->, bend left=55] (DA) to node {$\pi^0_n$}(A);
\end{tikzpicture}
}
\caption{Diagram of the main argument}
\label{Full diagram}
\end{figure}

Next we inductively define sequences $(M^i_n: i, n<\omega)$, $(\pi^i_n: i, n<\omega)$, $(\Lambda_i: i\leq \omega)$, $(\tau^{i, i+1}_n: i, n<\omega)$, $(\nu_n: i<\omega)$, $(\nu_n': i<\omega)$, $(\eta_n: n<\omega)$, $(\k_i:i<\omega)$, $(\theta_i: i<\omega)$, $(\T_i, E_i: i<\omega)$, $(M_i': i<\omega)$, $(\U_i, F_i: i<\omega)$, $(\P_i: i\leq n)$, $(\P_i': i<\omega)$, and $(\sigma_i: i<\omega)$ satisfying the following conditions (see Diagram \ref{Full diagram}).
\begin{enumerate}[(a)]
\item For all $i, n<\omega$, $\pi^i_n: M^i_n\rightarrow \P_i$ and $rng(\pi^i_n)\subseteq rng(\pi^i_{n+1})$.
\item $\tau^{i, i+1}_n: M^i_n\rightarrow M^{i+1}_n$. Let $\tau_n: M^0_n\rightarrow M^n_n$ be the composition of $\tau^{j, j+1}_n$'s for $j<n$. 
\item For all $i, n<\omega$, $\k_n=\tau_n(\k_0)$, $\eta_n=\tau_n(\eta_{0})$, $\nu_n=\tau_n(\nu_0)$ and $\nu_n'=\tau_n(\nu_0')$. 
\item For all $n<\omega$, $\T_n$ is an iteration of $M^n_n|\nu_n'$ above $\nu_n$ that makes $w_n$ generic and $M_n'$ is its last model. 
\item $\theta_n=\pi^{\T_n}(\nu_n')$ and $E_n\in \vec{E}^{M_n'}$ is such that $lh(E_n)>\theta_n$ and $\cp(E_n)=\k_n$.
\item for all $m,n$, $M^{n+1}_m=Ult(M^n_m, E_n)$ and $\tau_m^{n, n+1}=\pi_{E_n}^{M^n_m}$.
\item $\U_n=\pi^n_n\T_n$, $\P_n'$ is the last model of $\U_n$, $\sigma_n: M_n'\rightarrow \P_n'$ is the copy map and $F_n=\sigma_n(E_n)$.\footnote{So $\oplus_{i\leq n}\T_i$ and $\oplus_{i\leq n} \U_i$ are sealed iterations based on $\kappa$.}
\item $\P_{n+1}=Ult(\P_n, F_n)$ and $\psi^{n+1}_m:M^{n+1}_m\rightarrow \P_{n+1}$ is given by $\pi^{n+1}_m(\pi_{E_n}^{M^n_m}(f)(a))=\pi_{F_n}^{\P_n}(\pi^n_m(f))(\sigma_n(a))$.
\item $\Lambda_n=(\pi^n_n$-pullback of $(\Psi^{g*h}_\l)_{\P_n|\psi_n(\nu_n)})_{\eta_n, \nu_n}=(\sigma_n$-pullback of $(\Psi^{g*h}_\l)_{\P_n'|\sigma_n(\nu_n)})_{\eta_n, \nu_n}$ (see \rcor{coherence}).
\end{enumerate}

Let $M^{\omega}_{n}$ be the direct limit of $(M^m_n: m<\omega)$ under the maps $\tau^{m, m+1}_n$. Letting $\P_\omega$ be the direct limit of $(\P_n: n<\omega)$ and the compositions of $\pi_{F_n}^{\P_n}$, we have natural maps $\pi^\omega_n:M^\omega_n\rightarrow \P_\omega$. Notice that\\\\
(1) for each $n<\omega$, $\k_n<\omega_1^{V[g*h]}$ and $sup_n\k_n=\omega_1^{V[g*h]}$.\\\\
It follows that if $\tau^m_n: M_n^m\rightarrow M^\omega_n$ is the direct limit embedding then\\\\
(2) $\tau^m_n(\k_n)=\omega_1^{V[g*h]}$. \\\\
Next, notice that\\\\
(3) for each $m, n, p$, letting $\iota_n=\tau_n(\iota_{X_0})=\tau_n(\iota)$, $M^n_m|\iota_n=M^n_p|\iota_n$ and $\iota_n=(\k_n^+)^{M^n_m}$.\\
(4) for each $m,n, p$, $\pi^n_m\rest (M^n_m|\iota_n)=\pi^n_p\rest (M^n_p|\iota_n)$\\
(5) for each $m$, $n>1$ and $p>n$, $M^n_m|\theta_{n-1}=M^p_m|\theta_{n-1}$.\\
(6) for each $m$, $n>1$ and $p$ with $p>n$, $\pi^n_m\rest (M^n_m|\theta_{n-1})=\pi^p_m\rest (M^p_m|\theta_{n-1})$.\\\\
Because of condition (d) above we can find $G\subseteq Coll(\omega, <\omega_1^{V[g*h]})$ generic over $M^\omega_0$ such that $\bR^{M^\omega_0[G]}=\bR_{g*h}$ and $G\in V[g*h*k]$. By constructions, $\omega_1^{V[g*h]}$ is a limit of Woodin cardinals in $M^\omega_0$. It then follows from the results of \rsec{sec: gen-it} and \rsec{sec: cap ub} that 

\begin{lemma}\label{der model rep} $DM(G)^{M^\omega_0[G]}=L(\Gamma^\infty_{g*h}, \bR_{g*h})$. 
\end{lemma}
\begin{proof} It follows from \rcor{coherence} and \rlem{ub to strategy} that $A_n$ is projective in $\Lambda_n$. It follows from Corollary \ref{generic interpretability} that $\Lambda_n\rest HC^{V[g*h]} \in M^\omega_0[G]$ and it follows from \rcor{strategy in der model} that $\Lambda_n\rest HC^{V[g*h]}\in DM(G)^{M^\omega_0[G]}$. It follows that $\Gamma^\infty_{g*h}\subseteq DM(G)^{M^\omega_0[G]}$. 

 Moreover, it follows from \rcor{ub to strategy} that any set in $DM(G)^{M^\omega_0[G]}$ is projective in some $\Lambda_n\rest HC^{V[g*h]}$ and it follows from \rthm{strategies can be extended} that $\Lambda_n\rest HC^{V[g*h]}\in \Gamma^\infty_{g*h}$. Thus, $DM(G)^{M^\omega_0[G]}\subseteq L(\Gamma^\infty_{g*h}, \bR_{g*h})$. 
\end{proof}

We can also show variations of the above lemma for $M^\omega_n$ for each $n<\omega$. \rlem{der model rep} implies that in order to prove that $\sf{Sealing}$ holds, it is enough to establish clause 2 of $\sf{Sealing}$ as clause 1  immediately follows from \rlem{der model rep} and standard results about derived models (see \cite{DMT}). 

To continue, it will be easier to introduce some terminology. We say that the sequence $(X_i: i<\omega)$ is cofinal in $\Gamma^\infty_{g*h}$ as witnessed by $(A_i: i\in \omega)$ and $(w_i: i<\omega)$. We also say that $(M^n_0, \Lambda_n, \theta_n, \tau_{n, m}: n<m<\omega)$ is a $\Gamma^\infty_{g*h}$-genericity iteration induced by $(X_i: i<\omega)$ where $\tau_{n, m}: M^n_0\rightarrow M^m_0$ is the composition of $\tau^{i, i+1}_0$ for $i\in [n, m)$.  

\section{A proof of \rthm{main theorem}}\label{sec:main}

We now put together the results of the previous sections to obtain a proof of \rthm{main theorem}. Fix $h$ and $h'$ such that $h$ is $V[g]$-generic and $h'$ is $V[g*h]$-generic. We have shown in $V[g]$, clause (1) of $\sf{Sealing}$ holds. We now show clause (2) of $\sf{Sealing}$ holds in $V[g]$. We want to show that there is an embedding 
\begin{center}
$j: L(\Gamma_{g*h}, \bR_{g*h})\rightarrow L(\Gamma_{g*h*h'}, \bR_{g*h*h'})$
\end{center}
such that for $A\in \Gamma_{g*h}$, $j(A)=A^{h'}$. Let $(\xi_i: i<\omega)$ be an increasing sequence of cardinals such that $g*h*h'$ is generic for a poset in $V_{\xi_0}$. Let $u_n=(\xi_i: i< n)$. Set $W=L(\Gamma_{g*h}, \bR_{g*h})$ and $W'=L(\Gamma_{g*h*h'}, \bR_{g*h*h'})$.

Because $(\Gamma_{g*h})^{\#}$ exists, there is only one possibility for $j$ as above. Namely, given a term $\tau$, $n\in \omega$, $x\in \bR_{g*h}$ and $A\in \Gamma^\infty_{g*h}$, we must have that 
\begin{center}
$j(\tau^W(u_n, A, x))=\tau^{W'}(u_n, A^{h'}, x)$. 
\end{center}
What we must show is that $j$ is elementary. The next lemma finishes the proof. 

\begin{lemma}\label{j elem} $j$ is elementary.
\end{lemma}
\begin{proof} Let $u=(\eta, \d, \d', \l)$ be a good quadruple such that $\sup_{i<\omega}\xi_i<\eta$. Let $k\subseteq Coll(\omega, \Gamma^\infty_{g*h})$ be $V[g*h]$-generic and $k'\subseteq Coll(\omega, \Gamma^\infty_{g*h*h'})$ be $V[g*h*h']$-generic.

We have that $\Gamma^\infty_{g*h}$ is the Wadge closure of strategies of the countable substructures of $W_\l$. More precisely, given $A\in \Gamma^\infty_{g*h}$, there is an $X\prec (W_\l|\Psi_{\eta, \d}^{g*h})$ such that $A$ is Wadge reducible to $\Lambda_X$. It follows that to show that $j$ is elementary it is enough to show that given a formula $\phi$, $m\in \omega$, $X\prec ((W_\l, u)|\Psi_{\eta, \d}^{g*h})$ and a real $x\in \bR_{g*h}$,
\begin{center}
$W\models \phi[u_m, \Lambda_X, x]\Rightarrow W'\models \phi[u_m, \Lambda^{h'}_X, x]$.\footnote{The $\Leftarrow$ is similar as will be evident by the following proof.}
\end{center}
Fix then a tuple $(\phi, n, X, x)$ as above.

Working inside $V[g*h*k]$, let $(Y_i: i<\omega)$ be a cofinal sequence in $\Gamma^\infty_{g*h}$ as witnessed by some $\vec{A}$ and $\vec{w}$ such that $A_0=\emptyset$, $w_0=x$ and $Y'_0=X$. 

Working inside $V[g*h*h'*k']$, let $(Z_i: i<\omega)$ be a cofinal sequence in $\Gamma^\infty_{g*h*h'}$ as witnessed by some $\vec{B}$ and $\vec{v}$ such that $B_0=\emptyset$, $v_0=x$ and $Z'_0=X$. 

Let $(M_n, \Lambda_n, \theta_n, \tau_{n, l}: n<l<\omega)$ be a $\Gamma^\infty_{g*h}$-genericity iteration induced by $(Y_i: i<\omega)$ and $(N_n, \Phi_n, \nu_n, \sigma_{n, l}: n<l<\omega)$ be a $\Gamma^\infty_{g*h*h'}$-genericity iteration induced by $(Z_i: i<\omega)$. It is not hard to see that we can make sure that $M_1=N_1$ by simply selecting the same extender $E_0$ after $\T_0$; by our assumptions, $M_0 = N_0$ and $w_0 = v_0$. 

Let $\zeta=\eta_X$ and $\Gamma=(\Psi_{\eta, \d})_X$. Let $M_\omega$ be the direct limit along $(M_n: n<\omega)$ and $N_\omega$ the direct limit along $(N_n: n<\omega)$. For $n<\omega$, let $\k_n$ be the least strong cardinal of $M_n$ and $\k_n'$ be the least strong cardinal of $N_n$. Let $s_m^n$ be the first $m$ (cardinal) indiscernibles of $(M_n|\k_n)$ and $t_m^n$ be the first $m$ (cardinal) indiscernibles of $(N_n|\k_n')$. Notice that $(M_n|\k_n)^\#\in M_n$ and $(N_n|\k_n')^\#\in N_n$. It follows that 
$\tau_{n, l}(s^n_m)=s^l_m$ and $\sigma_{n, l}(t^n_m)=t^l_m$ for $n<l\leq \omega$. 

We then have the following sequence of implications. Below we let $\Gamma^*$ be the name for the generic extension of $\Gamma$ in the relevant model and $\dot{DM}$ be the name for the derived model. The third implication below uses the fact that $M_1 = N_1$.
\begin{align*}
W\models \phi[u_m, \Lambda_X, x] & \Rightarrow M_\omega[x]\models \emptyset \forces_{Coll(\omega, <\kappa_\omega)}\dot{DM}\models \phi[s^\omega_m, \Gamma^*, x]\\ &\Rightarrow M_1[x]\models \emptyset \forces_{Coll(\omega, <\kappa_1)}\dot{DM}\models \phi[s^1_m, \Gamma^*, x] \\ &\Rightarrow N_\omega[x]\models \emptyset \forces_{Coll(\omega, <\kappa'_\omega)}\dot{DM}\models \phi[t^\omega_m, \Gamma^*, x] \\ &\Rightarrow W'\models \phi[u_m, \Lambda^{h'}_X, x].
\end{align*}

\end{proof}

\section{$\sf{LSA-over-UB}$ may fail} \label{sec:LSAoverUBfails}

In this section, we prove Theorem \ref{thm:not_equiv}. We assume the hypotheses of \rthm{thm:not_equiv}. Here is the main consequence of the hypotheses that we need (see Lemma \ref{ub to strategy}):
\begin{enumerate}[(i)]
\item \label{one} letting $\lambda$ be an inaccessible cardinal which is a  limit of Woodin cardinals, and $g\subseteq Coll(\omega, <\lambda)$ be $V$-generic, then for any set $A$ which is Suslin co-Suslin in the derived model given by $g$, $DM(g)$ (see Section \ref{sec: cap ub}), then $A$ is Wadge reducible to $\Psi^g_{\delta}\rest HC^{V(\mathbb{R}^*_g)}$, for some Woodin cardinal $\delta < \lambda$. Furthermore, $\Psi^g_{\delta}\rest HC^{V(\mathbb{R}^*_g)}\in DM(g)$; in fact, $\Psi^g_{\delta}\rest HC^{V(\mathbb{R}^*_g)}=\Psi^g_{\delta}\rest HC^{V[g]}\in \Gamma^\infty_g$.
\end{enumerate}

 Suppose for contradiction that $\sf{LSA-over-UB}$ holds. Let $\lambda$ be an inaccessible cardinal which is a limit of Woodin cardinals in $V$. Let $h\subseteq Coll(\omega, <\lambda)$ be $V$ generic. By our assumption, in $V[h]$, there is some set $A$ such that
\begin{itemize}
\item $A\in V(\mathbb{R}^{V[h]})$;
\item $L(A,\mathbb{R})\models \sf{LSA}$;
\item $\Gamma^\infty_{h}$ is the Suslin co-Suslin sets of $L(A,\mathbb{R})$.
\item $\Gamma^\infty_{h} = \Delta_{h}$, where $\Delta_{h}$ is defined at the beginning of Section \ref{sec: cap ub}.
\end{itemize}
We note that the last item follows from (\ref{one}). 

Recall the notion of lbr hod mice is defined in \cite{normalization_comparison}. We will not need the precise definition of these objects. However, we need some notions related to short-tree strategies. Let $\P$ be a premouse (or hod premouse), $\tau$ a cut point cardinal of $\P$ (typically the $\tau$ we consider will be a Woodin cardinal or a limit of Woodin cardinals of $\P$), and $\Sigma$ an iteration strategy of $\P$ acting on trees based on $\P|\tau$. Suppose $\T$ according to $\Sigma$ is of successor length $\xi + 1$. Then we say $\T$ is \textit{short} if either $[0,\xi]_\T$ drops in model or else, letting $i$ be the branch embedding, $i(\tau) > \delta(\T)$; otherwise, we say $\T$ is \textit{maximal}. We let $\Sigma^{sh}$ be the short part of $\Sigma$; so $\Sigma^{sh}$ is a partial strategy. In the following, we may not have a (total) iteration strategy, but a partial strategy $\Lambda$ such that whenever $\T$ is according to $\Lambda$, if $\Lambda(\T)$ is defined, then either $\Lambda(\T)$ drops in model or else the branch embedding $i_{\Lambda(\T)}(\tau) > \delta(\T)$. We call such a $\Lambda$ a \textit{short-tree strategy}.\footnote{An example of a short-tree strategy is $\Sigma^{sh}$ for some total strategy $\Sigma$.} We may turn $\Lambda$ into a total strategy by assigning $\Lambda(\T)$ to be $\M(\T)^\sharp$ whenever a branch of $\T$ is not defined by $\Lambda$. Short tree strategies may be defined on stacks of normal trees as usual. 

The proof of \cite[Theorem 0.5]{HPC} gives us a pair $(\P,\Sigma)$ such that the following hold in $V(\mathbb{R}^{V[h]})$ (here the hypothesis $\sf{HPC+NLE}$ is applied in the model $L(A,\mathbb{R})$):
\begin{enumerate}
\item $\P$ is a least-branch hod premouse (lpm) (cf. \cite[Section 5]{normalization_comparison});
\item $\P$ has a largest Woodin cardinal $\delta = \delta^\P$ and letting $\kappa^\P$ be the least $<\delta$-strong cardinal in $\P$, then $\kappa^\P$ is a limit of Woodin cardinals;
\item \label{three}$\Sigma$ is a short-tree strategy of $\P$ and $\Sigma\in L(A,\mathbb{R})\backslash \Gamma^\infty_{h}$; furthermore, $\Sigma$ is Suslin in $L(A,\mathbb{R})$;
\item \label{four} for every $A\in \Gamma^\infty_{h}$, there is an iteration map $i:\P\rightarrow \Q$ according to $\Sigma$ such that $A <_w \Sigma_{Q|\kappa^\Q}$, where $\kappa^\Q$ is the least $\delta^\Q = i(\delta^\P)$-strong cardinal in $\Q$;
\item whenever $\T$ is according to $\Sigma$ and either $\Sigma(\T) = b$ is nondropping with last model $\Q$ or $\Sigma(\T)$ is not a branch with $\Q = \Sigma(\T)$, $\Sigma_{\T,\Q}$ satisfies (\ref{three}) and (\ref{four});\footnote{The above properties follow from the proof of Step 1 in \cite[Theorem 0.5]{HPC}, which can be applied to our hypothesis.}
\setcounter{nameOfYourChoice}{\value{enumi}}
\end{enumerate}
General properties of sets of reals in derived models give:
\begin{enumerate}
\setcounter{enumi}{\value{nameOfYourChoice}}
\item \label{six} there is some $\gamma < \lambda$ such that $(\P,\Sigma\rest V[h\rest \gamma]) \in V[h\rest \gamma]$.
\end{enumerate}

\begin{lemma}\label{lem:smaller}
Fix a $\gamma$ as in (\ref{six}). In $V[h\rest \gamma]$, there is a Woodin cardinal $\delta < \lambda$ such that $\delta > \gamma$ and there is a tree $\T$ according to $\Sigma$ such that either $\Sigma(\T)$ is a branch, $\Q = \M^\T_b$, and the branch embedding $i:\P\rightarrow \Q$ exists, or  $\Sigma(\T)$ is not a branch with $\Q = \Sigma(\T)$, and $\Sigma_{\T,\Q}$ satisfies (\ref{three}) and (\ref{four}) and is Wadge reducible to $\Psi^{h\rest \gamma}_\delta$.
\end{lemma}
\begin{proof}
Let $\delta$ be the least Woodin cardinal $> \gamma$. Let $\Psi = \Psi^{h\rest \gamma}_\delta$. Let $(\M_\xi, \Lambda_\xi: \xi \leq \delta)$ be the models and strategies of the fully backgrounded (lbr) hod mouse construction over $W^\Psi_\delta$ (cf. \cite{normalization_comparison}), where backgrounded extenders used have critical points $> max(\gamma, |\P|)$. Let $\T$ be according to $\Sigma$ be the comparison tree of $\P$ against the above construction. By universality, there is $\xi\leq \delta$ such that 
\begin{enumerate}[(i)]
\item either $\Sigma(\T) = b$ exists and there is an iteration map $i: \P \rightarrow \M_\xi$ and $\Sigma_{\T,\M_\xi} = \Lambda_\xi^{sh}$,
\item or $\Sigma(\T)$ does not exist ($\T$ is $\Sigma$-maximal), $\M_\xi = \Sigma(\T)$, and $\Sigma_{\T,\M_\xi} = \Lambda_\xi^{sh}$.
\end{enumerate}
In either case, we get that $\Sigma_{\T,\M_\xi}$ satisfies (\ref{three}) and (\ref{four}) above and $\Sigma_{\T,\M_\xi} = \Lambda_\xi^{sh}$ is Wadge reducible to $\Psi$ in $V(\mathbb{R}^{V[h]})$.\footnote{The fact that $\Lambda_\xi^{sh}$ is Wadge reducible to $\Psi$ is a standard property of fully backgrounded constructions. We abuse notations here, identifying for example $\Psi$ with its canonical extension in $V(\mathbb{R}^{V[h]})$.}
\end{proof}

Let $\delta, \T, \Q$ be as in Lemma \ref{lem:smaller}. Applying (\ref{one}) in $DM(h)$, we get that $\Psi^{h\rest \gamma}_\delta\rest HC^{V[h]}\in \Gamma^\infty_{h}$. Lemma \ref{lem:smaller} then implies that $\Sigma_{\T,\Q}\in \Gamma^\infty_{h}$. This contradicts (\ref{three}). This completes the proof of Theorem \ref{thm:not_equiv}.

\bibliographystyle{amsalpha}
\bibliography{SealingFromIter_submitted_revised.bib}
\end{document}